\documentclass{amsart}
\usepackage{times}
\usepackage{amsfonts,amssymb,amscd}
\usepackage{graphicx}

\newtheorem{thm}{Theorem}
\newtheorem{defi}{Definition}
\newtheorem{lem}{Lemma}
\newtheorem{cor}{Corollary}
\newtheorem{pro}{Proposition}
\newtheorem{rem}{Remark}

\newtheorem{obs}{Observation}
\newcommand{\pf}{{\bf\sl Proof.  }}
\renewcommand{\hbar}{\centerline{\rule{8cm}{0.5mm}}}
\renewcommand{\qed}{\vrule width0pt\hfill \raisebox{-.3ex}
   {\vrule height8pt width8pt depth0pt} \hspace*{-7pt}}
\renewcommand{\Re}{{\mathbb R}}
\newcommand{\C}{{\mathbb C}}

\newcommand{\tc}{{\mathcal C}_{\rm tot}}
\newcommand{\mc}{{\mathcal M}{\mathcal C}}

\newcommand{\nc}{{\mathcal N}}

\begin{document}

\large

\title[Total Curvature and Isotopy of Graphs]
{\textbf {Total Curvature and Isotopy of Graphs in $\Re^3$}}
\author{\bf Robert Gulliver and Sumio Yamada}
\thanks{Supported in part by JSPS Grant-in-aid for Scientific 
Research No.17740030}
\thanks{Thanks to the Korea Institute for Advanced
Study for invitations.}
\date{May 30, 2008}

\begin{abstract}
Knot theory is the study of isotopy classes of embeddings of the
circle $S^1$ into a $3$-manifold, specifically $\Re^3$.  The
F\'ary-Milnor Theorem says that any curve in $\Re^3$ of 
total curvature less
than $4\pi$ is unknotted.  More generally, a (finite) {\em graph}
consists of a finite number of edges and vertices.  Given a
topological type of graphs $\Gamma$, what limitations on the
isotopy class of $\Gamma$ are implied by a bound on total
curvature?  What does ``total curvature" {\em mean} for a graph?  We
define a natural notion of {\em net total curvature}
$\nc(\Gamma)$ of a graph $\Gamma \subset \Re^3$ (see
definition \ref{defnet} below), and prove
that if $\Gamma$ is homeomorphic to the $\theta$-graph, then
$\nc(\Gamma) \geq 3\pi$;  and if $\nc(\Gamma) < 4\pi$,
then $\Gamma$ is isotopic in $\Re^3$ to a planar $\theta$-graph.
Further, $\nc(\Gamma) = 3\pi$ only when $\Gamma$ is a convex
plane curve plus a chord.  We begin our discussion with piecewise
$C^2$ graphs, and extend all these results to continuous graphs in
the final section. In particular, we show that graphs of finite
total curvature are isotopic to polygonal graphs.
\end{abstract}

\maketitle

%
\section{Introduction:  Curvature of a Graph}\label{intro}

The celebrated F\'ary-Milnor theorem states that a curve in
$\Re^n$ of total curvature at most $4\pi$ is unknotted.  In the
present
paper, we shall emphasize the knotting dimension $n=3$.  As a key
step in his 1950 proof, John Milnor showed that for a smooth
Jordan curve $\Gamma$ in $\Re^3$, the total curvature equals
half the integral over $e \in S^2$ of the number $\mu(e)$ of
local maxima of the linear ``height" function $\langle e,\cdot
\rangle$
along $\Gamma$ \cite{M}.  This is then taken as the definition of
total curvature when $\Gamma$ is only continuous.
The F\'ary-Milnor theorem (even for $C^0$ curves)
follows, since total curvature less than $4\pi$ implies there is a
unit vector $e_0 \in S^2$ so that $\langle e_0,\cdot \rangle$ has
a unique local maximum, and therefore that this linear function
is increasing on an interval of $\Gamma$ and decreasing
on the complement.  Without changing the pointwise value of this
``height" function, $\Gamma$ can be topologically untwisted to a
standard embedding of $S^1$ into $\Re^3$.  The Fenchel theorem,
that any curve in $\Re^3$ has total curvature at least $2\pi$,
also follows from Milnor's key step, since for all $e\in S^2$, the
linear function $\langle e,\cdot \rangle$ assumes its maximum
somewhere along $\Gamma$, implying $\mu(e) \geq 1$.  Milnor's
proof is independent of the proof of Istvan F\'ary, published
earlier, which takes a different approach \cite{Fa}.

       We would like to extend these results, replacing the
simple closed curve by a finite {\em graph} $\Gamma$ in $\Re^3$.
$\Gamma$ consists of a finite number of points, called
{\em vertices}, and a finite number of simple arcs, called
{\em edges}, which each have as endpoints two of the vertices.  We
shall assume $\Gamma$ is connected.  The {\em valence} of a vertex
$q$ is the number $d(q)$ of edges which have $q$ as an endpoint.
For convenience, we shall assume in the first part of this paper
that each edge is of class $C^2$ up to its endpoints.  In section 
\ref{nonsmooth} below, we shall extend our main theorem to 
graphs which are merely {\it continuous}.  

\vspace{2mm}

\centerline{\bf{Definitions of Total Curvature for Knots}}

\vspace{2mm}

       The first difficulty, and a substantial one, in extending
the results of F\'ary and Milnor, is to understand the
contribution to total curvature at a vertex of valence
$d(q)\geq 3$.
For a smooth closed curve $\Gamma$, the total curvature is
$$
{\mathcal C}(\Gamma) = \int_\Gamma |\vec{k}| \, ds,
$$
where $s$ denotes arc length along $\Gamma$ and $\vec{k}$ is the
curvature vector.  If $x(s)\in \Re^3$ denotes the position of the
point measured at arc length $s$ along the curve, then
$\vec{k} = \frac {d^2x}{ds^2}$.  For a piecewise smooth curve,
that is, a graph with vertices $q_1, \dots, q_N$ having always
valence $d(q_i)=2$, the total curvature is readily generalized to
\begin{equation}\label{gencurv}
{\mathcal C}(\Gamma) =
\sum_{i=1}^N {\rm c}(q_i) + \int_{\Gamma_{\rm reg}} |\vec{k}| \, ds,
\end{equation}
where the integral is taken over the separate $C^2$ edges of
$\Gamma$ without their endpoints;
and where ${\rm c}(q_i) \in [0,\pi]$ is the
exterior angle formed by the two
edges of $\Gamma$ which meet at $q_i$.  That is,
$\cos({\rm c}(q_i)) = \langle T_1, -T_2\rangle,$
where $T_1= \frac{dx}{ds}(q_i^+)$ and
$T_2= -\frac{dx}{ds}(q_i^-)$ are the unit tangent vectors at
$q_i$ pointing into the two edges which meet at $q_i$.
The exterior angle ${\rm c}(q_i)$ is the correct contribution to total
curvature, since any sequence of smooth curves converging to
$\Gamma$ in $C^0$, with $C^1$ convergence on compact subsets of
each open edge, includes a small arc near $q_i$ along which the
tangent vector changes from near $\frac{dx}{ds}(q_i^-)$ to near
$\frac{dx}{ds}(q_i^+)$.  The greatest lower bound of the
contribution to total curvature of this disappearing
arc along the smooth approximating curves equals ${\rm c}(q_i)$.

In the context of the first variation of length of a graph,
however, a different notion of total curvature becomes
appropriate, known as the {\em mean curvature} of the graph
\cite{AA}.  The contribution to mean curvature at a vertex $q_i$ of
valence $2$, as above, is $|T_1 + T_2| = 2\sin(c(q_i)/2)$.  
This is smaller than the exterior angle $c(q_i)$.  If the
variation vector field is allowed to be multiple valued at $q_i$,
then the larger value $c(q_i)$ appears instead in the first 
variation of length.

\vspace{2mm}

\centerline{\bf{Definitions of Total Curvature for Graphs}}

\vspace{2mm}

When we turn our attention to a {\em graph} $\Gamma$, we find the
above definition for curves (valence $d(q)=2$) does not generalize
in any obvious way to higher valence (see \cite{G}).
The ambiguity of the general formula \eqref{gencurv} is resolved if
we specify the definition of ${\rm c}(0)$ when $\Gamma$
is the cone over a finite set $\{T_1, \dots, T_d\}$ in the unit
sphere $S^2$.  

A rather natural definition of total curvature of graphs was 
given by Taniyama in \cite{T}.  We have called this {\bf maximal 
total curvature} $\mc(\Gamma)$ in \cite{G};  the contribution to 
total curvature at a vertex $q$ of valence $d$ is
$${\rm mc}(q):=
\sum_{1\leq i<j\leq d}\arccos\langle T_i,-T_j\rangle.$$
In the case $d(q) = 2$, the sum above has only one term, the
exterior angle ${\rm c}(q)$ at $q$.
However, as may become apparent to readers of the present paper,
for $d(q) \geq 3$, this notion of total curvature 
offers certain problems for questions of the isotopy type 
of a graph in $\Re^3$.

In our earlier paper \cite{GY} on the density of an area-minimizing
two-dimensional 
rectifiable set $\Sigma$ spanning $\Gamma$, we found that it was
very useful to apply the Gauss-Bonnet formula to the cone over
$\Gamma$ with a point $p$ of \, $\Sigma$ as vertex.  The
relevant notion of total curvature in that context is {\bf cone
total curvature} $\tc(\Gamma)$, defined using ${\rm tc}(q)$ as the
choice for ${\rm c}(q)$ in equation \eqref{gencurv}:
\begin{equation}\label{deftc}
{\rm tc}(q) := \sup_{e \in S^2} \left\{
\sum_{i=1}^d\left(\frac{\pi}{2}-\arccos\langle T_i, e\rangle
\right) \right\}.
\end{equation}
Note that in the case $d(q) = 2$, the supremum above is assumed at
vectors $e$ lying in the smaller angle between the tangent vectors
$T_1$ and $T_2$ to $\Gamma$, so that ${\rm tc}(q)$ is then the
exterior angle ${\rm c}(q)$ at $q$.  The main result of \cite{GY} is
that
$2\pi$ times the area density of $\Sigma$ at any of its points is at
most equal to $\tc(\Gamma)$.  The same result had been proved by
Eckholm, White and Wienholtz for the case of a simple closed curve
\cite{EWW}.  Taking $\Sigma$ to be the branched immersion of the
disk given by Douglas \cite{D1} and Rad\'o \cite{R}, it follows
that if ${\mathcal C}(\Gamma) \leq 4\pi$, then $\Sigma$ is embedded, and
therefore $\Gamma$ is unknotted.  Thus \cite{EWW} provided an
independent proof of the F\'ary-Milnor theorem.  However,
$\tc(\Gamma)$ may be small for graphs which are far from the
simplest isotopy types of a graph $\Gamma$ (see the nearly
spherical example following Observation \ref{mc>>nc} below).

In this paper, we introduce the notion of {\bf net total
curvature}
${\mathcal N}(\Gamma)$, which is the appropriate definition for
generalizing --- {\em to graphs} --- Milnor's approach to isotopy
and total curvature of {\em curves}.  For each unit
tangent vector $T_i$ at $q$, $1 \leq i \leq d=d(q)$, let
$\chi_i:S^2 \rightarrow \{-1, +1\}$ be equal to $-1$ on the
hemisphere with center at $T_i$, and $+1$ on the opposite
hemisphere (values along the equator, which has measure zero, are
arbitrary).  We then define
%
%
\begin{equation}\label{defnc}
{\rm nc}(q):=
\frac{1}{4}\int_{S^2}\left[\sum_{i=1}^d\chi_i(e)\right]^+\,dA_{S^2}(e).
\end{equation}
In the case $d(q)=2$, the integrand is positive (and equals 2)
only on the set of unit vectors $e$ which have a negative inner
product with both $T_1$ and $T_2$, 
ignoring $e$ in sets of measure zero.
This set is a lune bounded by semi-great circles orthogonal to
$T_1$ and to $T_2$, and has spherical area equal to twice the
exterior angle.  So in this case, ${\rm nc}(q)$ is the exterior
angle.
%
%
\begin{defi}\label{defnet}
We define the {\em net total curvature} of a piecewise $C^2$
graph $\Gamma$ with vertices $\{q_1, \dots, q_N\}$ as
%
%
\begin{equation}
{\nc}(\Gamma):=
\sum_{i=1}^N {\rm nc}(q_i)+\int_{\Gamma_{\rm reg}} |\vec{k}| \, ds.
\end{equation}
\end{defi}

We would like to explain how the net total curvature
${\mathcal N}(\Gamma)$ of a graph is related to more familiar
notions of total curvature.  Recall that a graph $\Gamma$ has an
Euler circuit if and only if its vertices all have even valence,
by a theorem of Euler. 
An Euler circuit is a closed, connected path which traverses each
edge of
$\Gamma$ exactly once.  Of course, we do not have the hypothesis
of even valence.  We can attain that hypothesis by passing to the
{\em double} $\widetilde{\Gamma}$ of $\Gamma$:  let
$\widetilde{\Gamma}$ be the graph with the same vertices as
$\Gamma$, but with two copies of each edge of $\Gamma$.  Then at
each vertex $q$, the valence as a vertex of $\widetilde{\Gamma}$
is $\widetilde{d}(q) = 2\,d(q)$, which is even.  By Euler's theorem,
there is an Euler circuit $\Gamma'$ of \, $\widetilde{\Gamma}$, which
may be thought of as a closed path which traverses each edge of
$\Gamma$ exactly {\em twice}.  Now at each of the points
$\{q_1, \dots, q_d\}$ along $\Gamma'$ which are mapped to
$q \in \Gamma$, we may
consider the exterior angle ${\rm c}(q_i)$.  One-half the sum of these
exterior angles, however, depends on the choice of the Euler
circuit $\Gamma'$.  For example, if $\Gamma$ is the union of the
$x$-axis and the $y$-axis in Euclidean space $\Re^3$,
then one might choose $\Gamma'$ to have four right angles, or to
have four straight angles, or something in between, with
completely different values of total curvature.  In order to form
a version of total curvature at a vertex $q$ which only depends on
the original graph $\Gamma$ and not on the choice of Euler circuit
$\Gamma'$, it is necessary to consider some of the exterior
angles as partially balancing others.  In the example just
considered, where $\Gamma$ is the union of two orthogonal lines,
opposite right angles will be considered to balance each other
completely, so that ${\rm nc}(q)=0$, regardless of the choice of
Euler circuit of the double.

It will become apparent to the reader that the precise character
of an Euler circuit of $\widetilde\Gamma$ is not necessary in what
follows.  Instead, we shall refer to a {\em parameterization}
$\Gamma'$ of
the double $\widetilde\Gamma$, which is a mapping from a
$1$-dimensional manifold, not necessarily connected;  the mapping 
is assumed to cover each edge of $\widetilde\Gamma$ once.

The nature of ${\rm nc}(q)$ is clearer when it is localized on
$S^2$, analogously to
\cite{M}.  In the case $d(q)=2$, Milnor showed that the
exterior angle at the vertex $q$ equals $\frac{1}{2}$ the area of
those $e \in S^2$ such that the linear function
$\langle e, \cdot \rangle$, restricted to $\Gamma$, has a local
maximum at $q$.  In our context, we may describe ${\rm nc}(q)$ 
as one-half the integral over the sphere of the number of
{\em net local maxima}, which balances local maxima and local
minima against each other.  Along the parameterization $\Gamma'$ of
the double of $\Gamma$, the linear function
$\langle e, \cdot \rangle$ may have a local maximum at some of the
vertices $q_1, \dots, q_d$ over $q$, and may have a local minimum at
others.  In our construction, each local minimum balances against
one local maximum.  If there are more local minima than local
maxima, the number ${\rm nlm}(e,q)$, the net number of local
maxima, will be negative;  however, we will use only the positive
part $[{\rm nlm}(e,q)]^+$.

       We need to show that
$$ \int_{S^2} [{\rm nlm}(e,q)]^+ \,dA_{S^2}(e)$$
is independent of the choice of parameterization, and in fact is
equal to $2 \, {\rm nc}(q)$;  this will follow
from another way of computing ${\rm nlm}(e,q)$, in the next
section (see Corollary \ref{cor2} below).

%
\section{Some Combinatorics}\label{comb}
%
%

%
\begin{defi}\label{defnlm}
Let a parameterization \,$\Gamma'$ of the double of \, $\Gamma$ be
given.  Then a vertex $q$ of \,$\Gamma$ corresponds to a number of
vertices of \,$\Gamma'$, this number being exactly the valence
$d(q)$ of $q$ as a vertex of \,$\Gamma$.  If $q \in \Gamma$ is not
a vertex, then as the need arises, we may consider $q$ as a vertex
of valence $d(q) = 2$.  Let ${\rm lmax}(e,q)$ be the number of
local maxima of $\langle e, \cdot \rangle$ along \,$\Gamma'$ at
the points $q_1, \dots, q_d$ over $q$,  and similarly let 
${\rm lmin}(e,q)$ be the
number of local minima.  Finally, we define the number of
{\em net local maxima} of $\langle e, \cdot \rangle$ at $q$ to be
$${\rm nlm}(e,q) = \frac12[{\rm lmax}(e,q) - {\rm lmin}(e,q)]$$.
\end{defi}

%
\begin{rem}
The definition of ${\rm nlm}(e,q)$ appears to depend not only on
$\Gamma$ but on a choice of the
parameterization $\Gamma'$ of the double of
$\Gamma$;  but $\Gamma'$ is not unique, and indeed
${\rm lmax}(e,q)$ and ${\rm lmin}(e,q)$ may depend on the choice
of $\Gamma'$.  However, we shall see in Corollary \ref{cor1} 
below that the number of {\bf net} local maxima ${\rm nlm}(e,q)$
is in fact independent of $\Gamma'$.
\end{rem}

%
\begin{rem}
We have included the factor $\frac12$ in the definition of
${\rm nlm}(e,q)$ in order to agree with the difference of the
numbers of local maxima and minima along a parameterization 
of $\Gamma$ itself, in those cases where these numbers may be 
defined, namely, if $d(q)$ is even. 
\end{rem}

We shall {\bf assume} for the rest of this section that a unit
vector $e$ has been chosen, and that the linear ``height" function
$\langle e, \cdot \rangle$ has only a
finite number of critical points along $\Gamma$;  this excludes
$e$ belonging to a subset of $S^2$ of measure zero.
We shall also assume that the graph $\Gamma$ is subdivided to
include among the vertices all critical points of the linear
function $\langle e, \cdot \rangle$, with of course
valence $d(q) = 2$ if $q$ is an interior point of one of the
original edges of $\Gamma$.

%
\begin{defi}\label{updown}
Choose a unit vector $e$.  At a point $q \in \Gamma$ of valence
$d = d(q)$, let the {\em up-valence} $d^+ = d^+(e,q)$ be the
number of edges of \,$\Gamma$ with endpoint $q$ on which
$\langle e, \cdot \rangle$ is greater (``higher") than
$\langle e, q \rangle$, the ``height" of $q$.
Similarly, let the {\em down-valence} $d^-(e,q)$ be the number of
edges along which $\langle e, \cdot \rangle$ is less than its
value at $q$.  Note that $d(q) = d^+(e,q) + d^-(e,q)$, for almost all
$e$ in $S^2$.
\end{defi}

%
\begin{lem}\label{combin}
{\bf (Combinatorial Lemma)}
${\rm nlm}(e,q) = \frac12[d^-(e,q) - d^+(e,q)]$.
\end{lem}
\noindent
\pf
Let a parameterization $\Gamma'$ of the double of $\Gamma$ be
chosen, with respect to which ${\rm lmax}(e,q)$,
${\rm lmin}(e,q)$, ${\rm nlm}(e,q)$ and $d^\pm = d^\pm(e,q)$ are
defined.  Recall the assumption above, that $\Gamma$ has been
subdivided so that along each edge, the linear function
$\langle e, \cdot \rangle$ is strictly monotone.

Consider a vertex $q$ of $\Gamma$, of valence $d=d(q)$.
The parameterization $\Gamma'$ of the double of $\Gamma$ has
$2d$ edges with an endpoint among the points
$q_1, \dots, q_d$ which are mapped to $q \in \Gamma$.  There are
two copies of each edge of $\Gamma$, and on $2d^+$, resp. $2d^-$
of these, $\langle e, \cdot \rangle$ is greater resp. less than
$\langle e, q \rangle$.
But for each $1\leq i\leq d$, the parameterization
$\Gamma'$ has exactly two edges which meet at $q_i$.  Depending 
on the up/down character of the two edges of $\Gamma'$
which meet at $q_i$, $1\leq i\leq d$, we can count:\\
(+) If $\langle e, \cdot \rangle$ is greater than
$\langle e, q \rangle$ on both edges, then $q_i$ is a local
minimum point;  there are ${\rm lmin}(e,q)$ of these among
$q_1, \dots, q_d$.  \\
(-) If $\langle e, \cdot \rangle$ is less than
$\langle e, q \rangle$ on both edges, then $q_i$ is a local
maximum point;  there are ${\rm lmax}(e,q)$ of these.  \\
(0) In all remaining cases, the linear function
$\langle e, \cdot \rangle$ is greater than $\langle e, q \rangle$
along one edge and less along the other, in which case $q_i$ is
not
counted in computing ${\rm lmax}(e,q)$ nor ${\rm lmax}(e,q)$;
there are
$d(q)-{\rm lmax}(e,q)-{\rm lmin}(e,q)$ of these.

        Now count the individual edges of $\Gamma'$:  \\
(+) There are ${\ \rm lmin}(e,q)$ pairs of edges, each of which is
part of a local minimum, both of which are counted among the
$2 d^+(e,q)$ edges of $\Gamma'$ with
$\langle e, \cdot \rangle$ greater than $\langle e, q \rangle$.\\
(-) There are ${\ \rm lmax}(e,q)$ pairs of edges, each of which is
part of a local
maximum;  these are counted among the number $2d^-(e,q)$ of edges
of $\Gamma'$ with
$\langle e, \cdot \rangle$ less than $\langle e, q \rangle$.
Finally,\\
(0) there are $d(q)-{\rm lmax}(e,q)-{\rm lmin}(e,q)$ edges of
$\Gamma'$ which are not part of a local
maximum or minimum, with $\langle e, \cdot \rangle$
greater than $\langle e, q \rangle$;  and an equal number of edges
with $\langle e, \cdot \rangle$ less than $\langle e, q \rangle$.

Thus, the total number of these edges of $\Gamma'$ with
$\langle e, \cdot \rangle$ greater than $\langle e, q \rangle$ is
$$
2d^+=
2{\ \rm lmin}+(d-{\rm lmax}-{\rm lmin})=d+{\rm lmin}-{\rm lmax}.
$$
Similarly,
$$
2d^-=
2{\ \rm lmax}+(d-{\rm lmax}-{\rm lmin})=d+{\rm lmax}-{\rm lmin}.
$$
Subtracting gives the conclusion:
$$
{\rm nlm}(e,q):=
\frac{{\rm lmax}(e,q)-{\rm lmin}(e,q)}{2}=
\frac{d^-(e,q)-d^+(e,q)}{2}.
$$
\qed

%
\begin{cor}\label{cor1}
The number of net local maxima ${\rm nlm}(e,q)$ is independent of
the choice of parameterization $\Gamma'$ of the double of $\Gamma$.
\end{cor}
\pf
Given a direction $e\in S^2$, the up-valence and down-valence
$d^\pm(e,q)$ at a vertex $q\in \Gamma$ are defined
independently of the choice of $\Gamma'$.
\qed

%
\begin{cor}\label{cor2}
For any $q \in \Gamma$, we have
${\rm nc}(q) =
\frac12\int_{S^2} \Big[{\rm nlm}(e,q)\Big]^+ \,dA_{S^2}.$
\end{cor}
\pf
Consider $e \in S^2$.
In the definition \eqref{defnc} of ${\rm nc}(q),$ 
$\chi_i(e) = \pm 1$
whenever $\pm \langle e, T_i \rangle < 0$.  But the number of
$1\leq i \leq d$ with $\pm \langle e, T_i \rangle < 0$ equals
$d^{\mp}(e,q)$, so that
$$ \sum_{i=1}^d \chi_i(e) = d^-(e,q) - d^+(e,q)
= 2 \, {\rm nlm}(e,q)$$
by Lemma \ref{combin}.
\qed

%
\begin{defi}\label{defmu}
For a graph $\Gamma$ in $\Re^3$ and $e \in S^2$, define the
{\em multiplicity at $e$} as
$$ \mu(e) = \mu_\Gamma(e) = \sum\{{\rm nlm}^+(e,q): q
{\rm \ a\ vertex\ of\ } \Gamma
{\rm \ or\ a\ critical\ point\ of\ } \langle e,\cdot \rangle\}.$$
\end{defi}

Note that $\mu(e)$ is a half-integer.  Note also that in the case
when $\Gamma$ has no topological vertices, or
equivalently, when $d(q) \equiv 2$, $\mu(e)$ is exactly
the quantity $\mu(\Gamma, e)$, the number of local maxima of
$\langle e, \cdot \rangle$ along $\Gamma$ as defined in \cite{M},
p. 252.  It was shown in Theorem 3.1 of that paper that, in this
case, ${\mathcal C}(\Gamma) = \frac12 \int_{S^2} \mu(e) \,
dA_{S^2}$.
We may now extend this result to {\em graphs}:

%
\begin{thm}\label{muthm}
For a (piecewise $C^2$) graph $\Gamma,$ the net total curvature 
has the following representation:
$$ \nc(\Gamma) = \frac12 \int_{S^2} \mu(e) \,dA_{S^2}(e). $$
\end{thm}
\pf
We have $ \nc(\Gamma) =
\sum_{j=1}^N {\rm nc}(q_j) + 
\int_{\Gamma_{\rm reg}} |\vec{k}| \, ds,$
where $q_1, \dots, q_N$ are the vertices of \, $\Gamma$, including
vertices of valence $d(q_j) = 2$, and where
$\rm{nc}(q):=
\frac12\int_{S^2}\left[\sum_{i=1}^d\chi_i(e)\right]^+\,dA_{S^2}(e)$
by Definition \ref{defnc}.  By Milnor's result,
${\mathcal C}(\Gamma_{\rm reg}) =
\frac12 \int_{S^2} \mu_{\Gamma_{\rm reg}}(e) \, dA_{S^2}$.  But
$\mu_{\Gamma}(e) = \mu_{\Gamma_{\rm reg}}(e) +
\sum_{j=1}^N \rm{nlm}^+(e,q_j)$, and the theorem follows.
\qed\\

In section \ref{nonsmooth}, we will have need of the monotonicity 
of $\nc(P)$
under refinement of {\em polygonal} graphs $P$. This follows from 
the following Proposition.

%
\begin{pro}\label{monotmu}
Let $P$ and $\widetilde{P}$ be polygonal graphs in $\Re^3$, having
the same topological vertices, and homeomorphic to each other.
Suppose that every vertex of $P$ is also a vertex of $\widetilde{P}$.
Then for almost all $e \in S^2$, the multiplicity 
$\mu_{\widetilde{P}}(e) \geq \mu_P(e).$  As a consequence, 
$\nc(\widetilde{P}) \geq \nc(P)$.
\end{pro}
\pf
We may assume, as an induction step, that $\widetilde{P}$ 
is obtained from $P$ by replacing
the edge having endpoints $q_0$, $q_2$ with two edges, one having
endpoints $q_0$, $q_1$ and the other having endpoints $q_1$, $q_2$. 
Choose $e \in S^2$.  We consider several cases:

If the new vertex $q_1$ satisfies 
$\langle e, q_0\rangle < \langle e, q_1\rangle < 
\langle e, q_2\rangle$, then 
${\rm nlm}_{\widetilde{P}}(e,q_i)={\rm nlm}_P(e,q_i)$ for 
$i = 0,2$ and ${\rm nlm}_{\widetilde{P}}(e,q_1)=0$, hence
$\mu_{\widetilde{P}}(e) = \mu_P(e)$. 

If $\langle e, q_0\rangle < \langle e, q_2\rangle < 
\langle e, q_1\rangle$, then 
${\rm nlm}_{\widetilde{P}}(e,q_0)={\rm nlm}_P(e,q_0)$ and
${\rm nlm}_{\widetilde{P}}(e,q_1)=1$.  The vertex $q_2$ 
requires more careful counting:  the up- and down-valence 
$d_{\widetilde{P}}^\pm(e,q_2)=d_P^\pm(e,q_2) \pm 1$, so that by
Lemma \ref{combin}, 
${\rm nlm}_{\widetilde{P}}(e,q_2)={\rm nlm}_P(e,q_2)-1$.
Meanwhile, for each of the polygonal graphs, $\mu(e)$ is the sum
over $q$ of ${\rm nlm}^+(e,q)$, so the change from $\mu_P(e)$ to
$\mu_{\widetilde{P}}(e)$ depends on the value of 
${\rm nlm}_P(e,q_2)$:\\
(a) if ${\rm nlm}_P(e,q_2)\leq 0$, then 
${\rm nlm}_{\widetilde{P}}^+(e,q_2)={\rm nlm}_P^+(e,q_2)=0$;\\
(b) if ${\rm nlm}_P(e,q_2) = \frac12,$ then 
${\rm nlm}_{\widetilde{P}}^+(e,q_2)=
{\rm nlm}_P^+(e,q_2)-\frac12$;\\
(c) if ${\rm nlm}_P(e,q_2)\geq 1$, then 
${\rm nlm}_{\widetilde{P}}^+(e,q_2)=
{\rm nlm}_P^+(e,q_2)-1$.\\
Since the new vertex $q_1$ does not appear in $P$, recalling
that ${\rm nlm}_{\widetilde{P}}(e,q_1)=1$, we have 
$\mu_{\widetilde{P}}(e) - \mu_P(e) = +1, +\frac12$ or $0$ in the
respective cases (a), (b) or (c).  In any case, 
$\mu_{\widetilde{P}}(e) \geq \mu_P(e)$. 

The reverse inequality
$\langle e, q_1\rangle < \langle e, q_2\rangle < 
\langle e, q_0\rangle$ leads to a similar case-by-case argument. 
This time, ${\rm nlm}_{\widetilde{P}}(e, q_1) = -1$, so $q_1$ does
not contribute to $\mu_{\widetilde{P}}(e)$.  We have 
${\rm nlm}_{\widetilde{P}}(e, q_0)-{\rm nlm}_P(e,q_0) = 0$,
while 
${\rm nlm}_{\widetilde{P}}(e, q_2)-{\rm nlm}_P(e,q_2) =
\frac12[d_{\widetilde{P}}^--d_{\widetilde{P}}^+-d_P^-+ d_P^+]=1$.
Now depending whether ${\rm nlm}_P(e,q_2)$ is $\leq -1$,
$=-\frac12$ or $\geq 0$, we find that 
$\mu_{\widetilde{P}}(e)-\mu_P(e)= 
{\rm nlm}^+_{\widetilde{P}}(e, q_2)-{\rm nlm}^+_P(e,q_2) =
0$, $\frac12$, or $1$.  In any case,
$\mu_{\widetilde{P}}(e)\geq \mu_P(e)$.

Note that these arguments are unchanged if $q_0$ is switched with
$q_2$.  This covers all cases except those in which equality
occurs between $\langle e,q_i \rangle$ and $\langle e,q_j \rangle$
($i\neq j$).  The set of such unit vectors $e$ form a set of
measure zero in $S^2$.  The conclusion 
$\nc(\widetilde{P}) \geq \nc(P)$ now follows from Theorem
\ref{muthm}.  \qed 

%
%
\section{Valence three or four}\label{three/four}

Before proceeding to 
the relation between $\nc(\Gamma)$ and isotopy, we shall illustrate
some properties of net total curvature $\nc(\Gamma)$ in a 
few relatively simple cases.

%
\subsection{Minimum curvature for valence three}

%
\begin{pro}\label{val3}
If a vertex $q$ has valence $d(q)=3$, then 
$\rm{nc}(q) \geq \pi/2$,
with equality if and only if the three tangent vectors
$T_1, T_2, T_3$ at $q$ are coplanar but do not lie in any 
open half-plane.
\end{pro}
\pf At a vertex of valence $3$, we have the simplification that
the local parameterization $\Gamma'$ of the double
$\widetilde\Gamma$
is unique.  Write $\alpha_i \in [0,\pi]$ for the angle between
$T_i$ and $T_{i+1}$ (subscripts modulo $3$).  Then
$\rm{nc}(q) = \frac12 \sum_{i=1}^3 (\pi - \alpha_i)=
\frac{3}{2}\pi -  \frac12 \sum_{i=1}^3 \alpha_i$;  but
$\sum_{i=1}^3 \alpha_i \leq 2\pi,$ with strict inequality if the
vectors are not coplanar or lie in an open half-plane.  Otherwise,
equality holds.
\qed

It might be observed that among vertices of a given 
valence $d$, the minimum of $\rm{nc}(q)$ is $0$ if $d$ is even.

%
\subsection{Non-monotonicity of $\nc$ for subgraphs}

%
\begin{obs}\label{notmonotone}
It might be assumed that if $\Gamma_0$ is a subgraph of a graph
$\Gamma$, then $\nc(\Gamma_0) \leq \nc(\Gamma).$  However, this
is
{\bf not} the case.
\end{obs}

For a simple polyhedral counterexample, we may
consider the ``butterfly" graph $\Gamma$ in the plane with six
vertices:  $q_0^\pm = (0,\pm 1), q_1^\pm = (1,\pm 3),$ and
$q_2^\pm = (-1,\pm 3)$.  Three vertical edges $L_0, L_1$ and $L_2$
are the line segments $L_i$ joining $q_i^-$ to $q_i^+$;  four
additional edges are the line segments from $q_0^\pm$ to $q_1^\pm$
and from $q_0^\pm$ to $q_2^\pm$.  The small angle $2 \alpha$
at $q_0^\pm$ has $\tan \alpha = 1/2$, so that $\alpha < \pi/4.$

The subgraph $\Gamma_0$ will be $\Gamma$ minus the interior of
$L_0$.  Then $\Gamma_0$ is a simple closed curve, so that
at each vertex $q_i^\pm$, we have
$\rm{nc}_{\Gamma_0}(q) = {\rm c}(q) \in [0,\pi]$,
the exterior angle.  The edges are all straight, so the net total
curvature has contributions only at the six vertices;  we compute
$$\nc(\Gamma_0) = {\mathcal C}(\Gamma_0)=
4(\pi - \alpha) + 2(\pi - 2 \alpha) = 6 \pi - 8 \alpha.$$

Meanwhile, as we have seen in Proposition \ref{val3}, the three
(coplanar) edges of $\Gamma$ at each of the two vertices $q_0^\pm$
determine $\rm{nc}_\Gamma(q_0^\pm) = \pi/2$, so that
$$\nc(\Gamma) =
4(\pi - \alpha) + 2(\pi/2) = 5 \pi - 4 \alpha.$$

Since $\alpha < \pi/4,$ this implies that
$\nc(\Gamma_0) > \nc(\Gamma).$
\qed

Monotonicity in the sense of Observation \ref{notmonotone} is a
virtue of Taniyama's total curvature $\mc(\Gamma)$.

%
\subsection{Simple description of net total
curvature for valence $3$}

%
\begin{pro}\label{net3}
If $\, \Gamma$ is a graph having vertices only of valence two or
three, then $\nc(\Gamma) = \frac12 {\mathcal C}(\Gamma')$ for any
parameterization $\Gamma'$ of the double of $\, \Gamma$.
\end{pro}
\pf
Since $\Gamma'$ covers each edge of $\Gamma$ twice, we need only
show, for every vertex $q$ of $\Gamma$, having valence
$d(q) \in \{2,3\}$, that
%
%
\begin{equation}\label{*}
2\, {\rm nc}_\Gamma(q)= \sum_{i=1}^d {\rm c}_{\Gamma'}(q_i),
\end{equation}
where $q_1, \dots, q_d$ are the vertices of $\Gamma'$ over $q$.
If $d=2$, then
${\rm nc}_\Gamma(q)={\rm c}_{\Gamma'}(q_1)=
{\rm c}_{\Gamma'}(q_2)$, so equation \eqref{*} clearly holds.  For
$d=3$, write both sides of equation \eqref{*} as integrals over
$S^2$, using the definition \eqref{defnc} of ${\rm nc}_\Gamma(q)$:
we need to show that
\begin{eqnarray*}
2\int_{S^2}\left[\chi_1+\chi_2+\chi_3\right]^+\,dA_{S^2}&=&
\int_{S^2}\left[\chi_1+\chi_2\right]^+\,dA_{S^2}+ \\
+ \int_{S^2}\left[\chi_2+\chi_3\right]^+\,dA_{S^2}&+&
\int_{S^2}\left[\chi_3+\chi_1\right]^+\,dA_{S^2},
\end{eqnarray*}
where at each direction $e\in S^2$, $\chi_i(e) = \pm 1$ is the
sign of $\langle -e, T_i\rangle$.  But the integrands are equal
at almost every point $e$ of $S^2$:
$$
2\left[\chi_1+\chi_2+\chi_3\right]^+ =
\left[\chi_1+\chi_2\right]^+ +
\left[\chi_2+\chi_3\right]^+ +
\left[\chi_3+\chi_1\right]^+,
$$
as may be confirmed by cases:  $6=6$ if $\chi_1=\chi_2=\chi_3=+1$;
$2=2$ if exactly one of the $\chi_i$ equals $-1$ and $0=0$ in the
remaining cases.
\qed

%
\subsection{Simple description of net total
curvature fails, $d \geq 4$}

%
\begin{obs}\label{notinf}
We have seen, for any parameterization $\Gamma'$ of
the double $\widetilde\Gamma$ of a graph $\Gamma$, that
$\nc(\Gamma) \leq \frac12 {\mathcal C}(\Gamma')$, the total
curvature in the usual sense of the parameterized curve
$\Gamma'$.  Moreover, for graphs with vertices of valence $\leq 3$,
equality holds, by Proposition \ref{net3}.
A natural suggestion would be that $\nc(\Gamma)$
might be half the infimum of total curvature of all parameterizations
$\Gamma'$ of the double. However, in some cases, we have the
{\bf strict inequality}
$\nc(\Gamma) < \inf_{\Gamma'}\frac12 {\mathcal C}(\Gamma')$.
\end{obs}

In light of Proposition \ref{net3},
we choose an example of a vertex $q$ of valence four.

Suppose that for a small positive angle $\alpha$,
($\alpha \leq 1$ radian would suffice)
the four unit tangent vectors at $q$ are $T_1 = (1,0,0)$;
$T_2 = (0,1,0)$; $T_3 = (-\cos\alpha,0,\sin\alpha)$; and
$T_4 = (0,-\cos\alpha,-\sin\alpha)$.
Write $\theta_{ij} = \pi - \arccos \langle T_i, T_j \rangle.$
Then each of the 
possible parameterizations
of the double $\widetilde\Gamma$ has total curvature
equal to the sum of any four of the $\theta_{ij}$, where
each of the subscripts $1,2,3$ and $4$ appears twice.
The particular parameterization where the four connected
components are
criss-crossing at the origin, which has total exterior curvature
$2\theta_{13}+2\theta_{24} = 4\alpha$,
realizes the minimum value of the exterior curvature,
as it can be seen as a perturbation of the simple ``X-crossing"
case ($\alpha=0$) in which case there is
no curvature contribution at the vertex.  This shows that
$\inf_{\Gamma'}\frac12 {\mathcal C}(\Gamma') = 2\alpha.$

However, ${\rm nc}(q)$ is strictly less than $2\alpha$.  We have
written it as an integral over the unit sphere $S^2$:
$${\rm nc}(q):=
\frac12\int_{S^2}\left[\sum_{i=1}^d\chi_i(e)\right]^+\,dA_{S^2}(e)$$
(see equation \eqref{defnc}).
Note that $\chi_1(e)$ cancels $\chi_3(e)$, and $\chi_2(e)$ cancels
$\chi_4(e)$, except on regions bounded by two pairs of great circles:
$T_1^\perp$ and $T_3^\perp$;  and by $T_2^\perp$ and $T_4^\perp$.
Each of the four lune-shaped sectors has area $2\alpha$ (see
Figure~1).

\begin{figure}[ht]
\begin{center}
\includegraphics[scale=0.3]{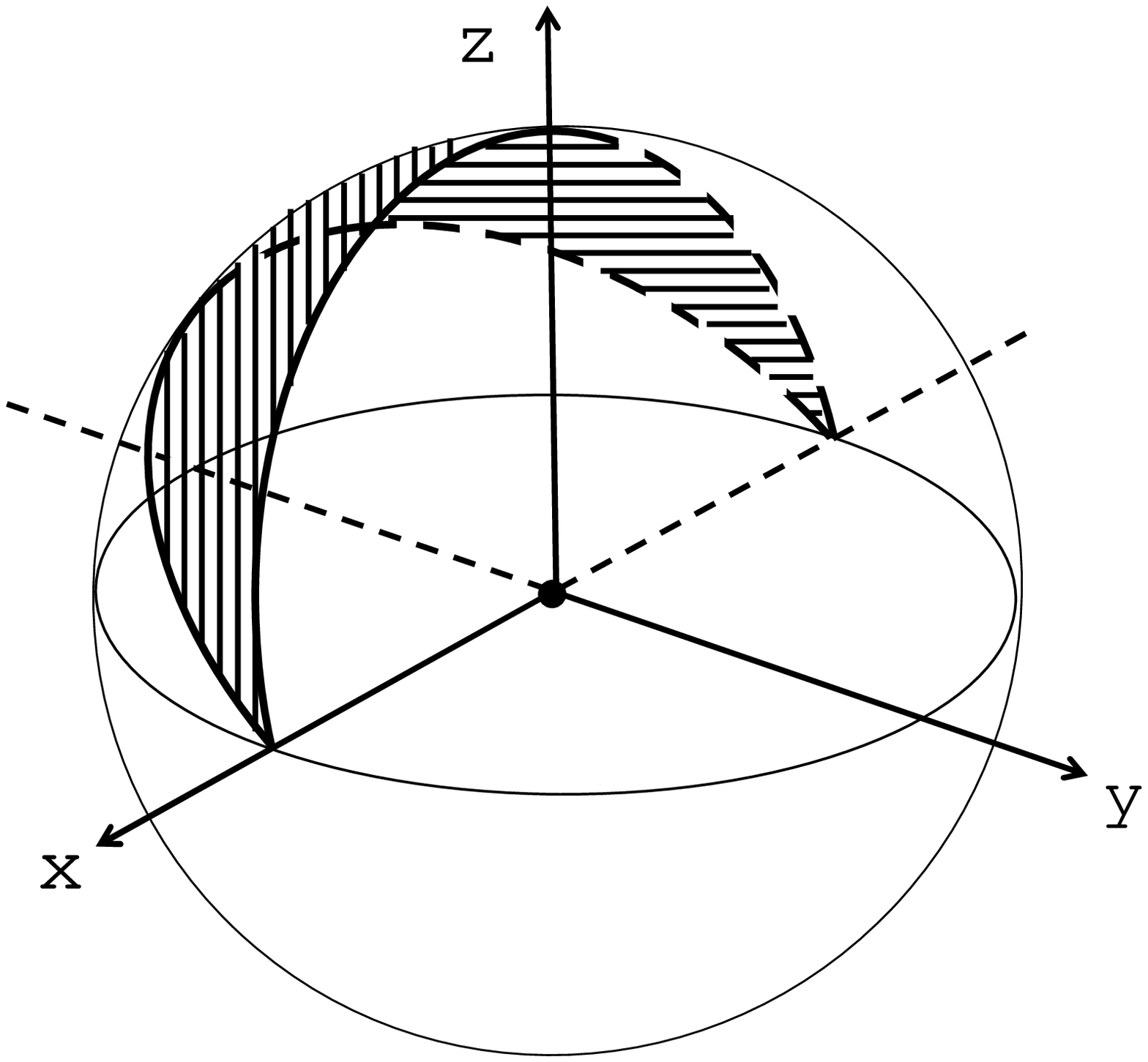}\quad
\includegraphics[scale=0.3]{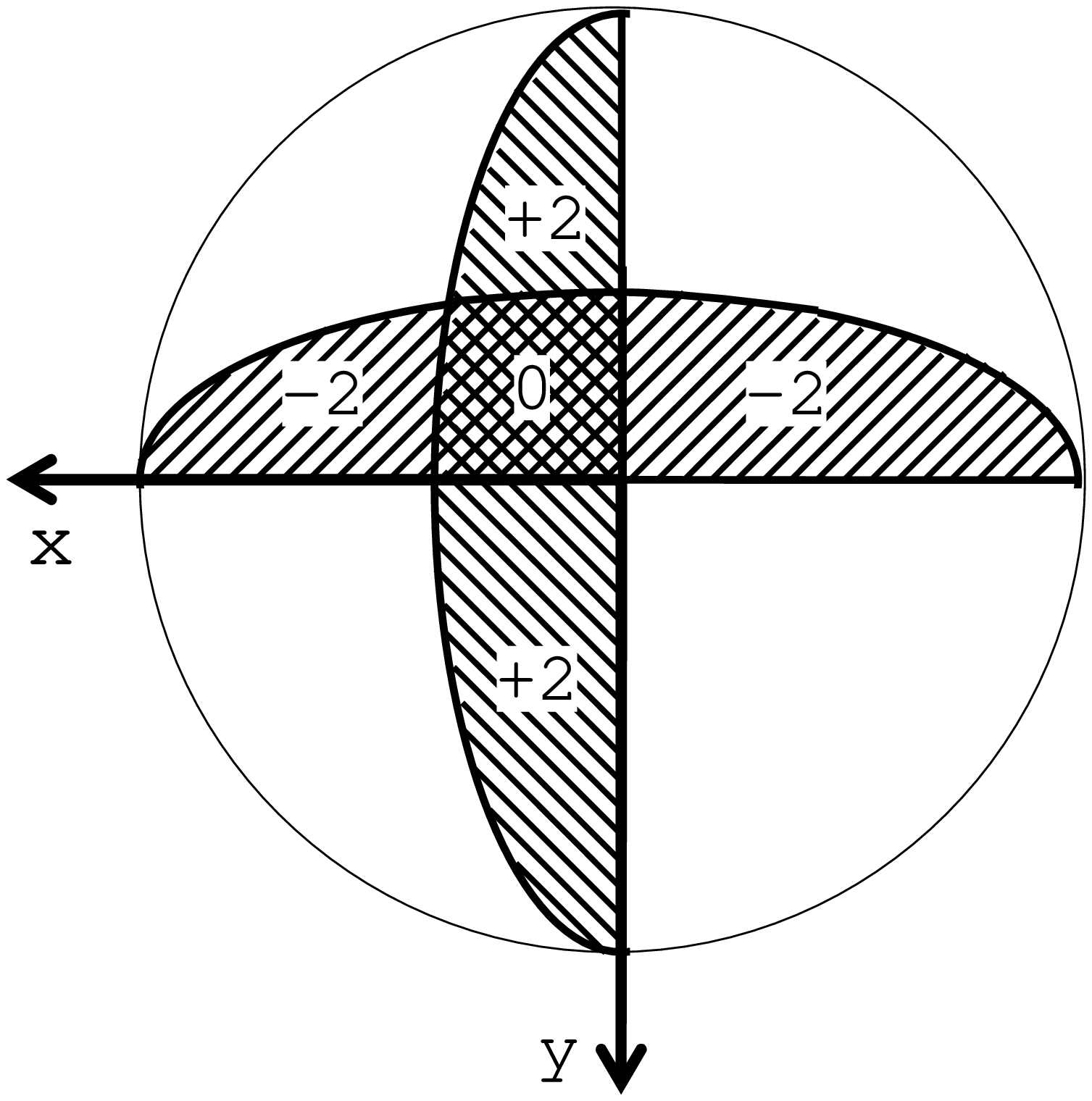}
\end{center}
\caption{an individual lune and overlapping lunes from above}
\label{figure1&2}
\end{figure}

Inside the sector partly bounded by the meridian
$T_1^\perp \cap \{z < 0\}$, the sum $\chi_1 + \chi_3 = 2,$
with the opposite sign on the antipodal sector.  The sum
$\chi_2 + \chi_4$ has a congruent
pattern, but rotated by a right angle so that the
sector where $\chi_1 + \chi_3$ is {\it positive} is
crossed by the sector where $\chi_2 + \chi_4$ is {\it negative}
(see Figure~1).
Thus, since we compute the integral only of the positive part of
the sum, there is partial cancellation $\sim \alpha^2$ when
$\alpha>0$ is small. It follows that 
$ {\rm nc}(q) < 2 \alpha $. \qed

%
\subsection{Net total curvature $\neq$ cone total curvature
$\neq$ Taniyama's total curvature}

It is not difficult to construct three unit
vectors $T_1, T_2, T_3$ such that the values of ${\rm nc}(q)$,
${\rm tc}(q)$ and ${\rm mc}(q)$, with these vectors as the 
$d(q) = 3$ tangent vectors to a graph at a vertex $q$, have quite 
different values.  For example, we
may take $T_1, T_2$ and $T_3$ to be
three unit vectors in a plane, making equal angles $2\pi/3$.
According to Proposition \ref{val3}, we have the contribution to
net total curvature ${\rm nc}(q) = \pi/2$.  But the contribution
to cone total curvature is ${\rm tc}(q) = 0$. Namely,
${\rm tc}(q) := \sup_{e \in S^2}
\sum_{i=1}^3\left(\frac{\pi}{2}-\arccos\langle T_i, e\rangle
\right).$  In this supremum, we may choose $e$ to be normal 
to the plane of $T_1, T_2$ and $T_3$, and ${\rm tc}(q) = 0$
follows.
Meanwhile, ${\rm mc}(q)$ is the sum of the exterior angles formed
by the three pairs of vectors, each equal to $\pi/3$, so that
${\rm mc}(q) = \pi$.

A similar computation for valence $d$ and coplanar vectors making
equal angles gives ${\rm tc}(q) = 0$, and 
${\rm mc}(q) = \frac{\pi}{2}\Big[\frac{(d-1)^2}{2}\Big]$ (brackets
denoting integer part), while ${\rm nc}(q) = \pi/2$ for $d$ odd,
${\rm nc}(q) = 0$ for $d$ even.  This example indicates that 
${\rm mc}(q)$ may be significantly larger than ${\rm nc}(q)$.
In fact, we have 

%
\begin{obs}\label{mc>>nc}
If a vertex $q$ of a graph $\Gamma$ has valence $d=d(q)$, then
${\rm mc}(q) \geq (d-1) {\rm nc}(q)$.
\end{obs}

This observation follows from the definition \eqref{defnc} of 
${\rm nc}(q)$, and the remark that if $d=2$, then ${\rm nc}(q)$
and ${\rm mc(q)}$ are both equal to
the exterior angle at $q$.  The case $d=3$ has
essentially already been treated in Proposition \ref{net3}.  In
the general case, let 
$T_1, \dots, T_d$ be the unit tangent vectors at $q$.  The 
exterior angle between $T_i$ and $T_j$ is 
$$\arccos\langle -T_i,T_j \rangle =
\frac{1}{4} \int_{S^2} (\chi_i + \chi_j)^+ \, dA_{S^2}.$$
The contribution ${\rm mc}(q)$ at $q$ to Taniyama's total curvature 
$\mc(\Gamma)$ equals the sum
of these integrals over all $1 \leq i < j \leq d$.  The sum of the
integrands is 
$$ \sum_{1\leq i<j\leq d}(\chi_i + \chi_j)^+\geq
\Big[\sum_{1\leq i<j\leq d}(\chi_i + \chi_j)\Big]^+ = 
(d-1)\Big[\sum_{i=1}^d \chi_i\Big]^+.$$
Integrating over $S^2$ and dividing by $4$, we have 
${\rm mc}(q) \geq (d-1) {\rm nc}(q)$.
\qed
\\

We may apply these computations
for a vertex with coplanar tangent vectors  making equal angles
to compare the three versions of curvature discussed in Section
\ref{intro} above, for a 
graph $\Gamma$ in $\Re^3$ in a specified isotopy class.  Any isotopy 
type may be represented in the sphere $S^2$ (or the plane) by a
piecewise smooth diagram in which at most two points are mapped to
the same point of $S^2$, the two points being interior points of
edges of $\Gamma$, with the over/under character 
of each crossing indicated.  At each vertex, 
let the diagram be given on the sphere with equal angles.  We
interpret the diagram in terms of a smooth function on the graph
$\Gamma$, with greater resp. smaller values along an edge passing
over resp.  under another edge.  Multiply this function by a small
positive constant $\varepsilon$ and add the constant $1$, to give
a representative in spherical coordinates of the given isotopy
class.  As $\varepsilon \to 0$, the contribution 
to cone total curvature at each vertex tends to $0$, but any
vertex of odd valence $d$ contributes $\frac{\pi}{2}$ in the limit
to net total curvature, while ${\rm mc}(q)$ is roughly
$\frac{\pi}{4}d^2$.  For a graph $\Gamma$ of a complicated isotopy
type, $\tc(\Gamma)$ will be approximately as large as the total
curvature of a similarly complicated knot;  for $\Gamma$ with only
vertices of even valence, $\nc(\Gamma)$ will equal $\tc(\Gamma)$ 
asymptotically; but where vertices
occur with high odd valence, $\nc(\Gamma)$ and especially
$\mc(\Gamma)$ may be much larger.  

%
%
\section{Isotopy of theta graphs}

We shall consider in this section one of the simpler homeomorphism
types of graphs, the {\bf theta graph.}  A theta graph has only
two vertices $q^\pm,$ and three edges, each of which
connects $q^+$ to $q^-$.  The {\bf standard theta graph} is the
isotopy class in $\Re^3$ of a plane circle plus a diameter, or of
the printed lower-case Greek letter $\theta$.  We shall show that
a theta graph of net total curvature $< 4\pi$ is isotopically
standard.

To simplify the exposition and highlight the role of net total
curvature, we first treat the case of piecewise $C^2$ graphs with
distinct tangent vectors at each vertex.  This theorem will be
generalized in Theorem \ref{theta2} below to the case  of arbitrary 
continuous theta graphs.

We may observe that there are nonstandard theta graphs in
$\Re^3$.  For example, the union of two edges might form a knot.

Using the notion of net total curvature, we may extend the
theorems of Fenchel \cite{Fen} as well as of F\'ary-Milnor
(\cite{Fa},\cite{M}), for curves homeomorpic to $S^1$, to graphs
homeomorphic to the theta graph.  We would like to thank Jaigyoung
Choe and Rob Kusner for their comments about the case
$\nc(\Gamma) = 3\pi$.

%
\begin{thm}\label{thetathm}
Suppose $\Gamma \subset \Re^3$ is a piecewise $C^2$ theta graph,
with distinct unit tangent vectors at each vertex.  
Then $\nc(\Gamma) \geq 3\pi$.  If $\nc(\Gamma) < 4\pi$, then
$\Gamma$ is isotopic in $\Re^3$ to the standard theta graph.
Moreover if $\nc(\Gamma) = 3\pi$, then the graph is a planar 
convex curve plus a straight chord.
\end{thm}

Recall Observation 1:  for a subgraph $\Gamma_0$ of $\Gamma$,
$\nc(\Gamma_0)$ may be greater than $\nc(\Gamma)$.  Thus, 
when $\nc(\Gamma) < 4\pi$, the
unknottedness of subgraphs needs proof.

\pf
We first show the Fenchel-type {\bf lower bound} $3\pi$.  For any
$e\in S^2$, consider the linear function $\langle \cdot,e \rangle$
along $\Gamma$.  We need to show that $\mu(e) \geq 3/2$, from
which Theorem \ref{muthm} implies that $\nc(\Gamma) \geq 3\pi$.
Write the two vertices of $\Gamma$ as $q^\pm$, each with valence
$d(q^\pm) = 3$.

If the maximum $t_{\rm max}$ along $\Gamma$ of
$\langle \cdot,e \rangle$ occurs at a vertex, say $q^+$, then
${\rm nlm}(e,q^+) = 3/2.$  Namely, at each of the three vertices
$q_1^+, q_2^+, q_3^+$ of $\Gamma'$ corresponding to $q^+$,
$\langle \cdot,e \rangle$ has a (global hence local) maximum.
It follows in this case from Definitions \ref{defnlm} and
\ref{defmu} that $\mu(e) \geq 3/2$.

Otherwise, the maximum occurs at an interior point $q_{\rm max}$
of one of the topological edges of $\Gamma$; then
${\rm nlm}(e,q_{\rm max}) = 1$.  By disregarding a set of
$e \in S^2$ of measure zero, we may assume $q_{\rm max}$ is the
unique maximum point.  Then for values of $t$ slightly smaller
than $t_{\rm max} = \langle q_{\rm max},e \rangle$, there
are exactly two points, which will be close to $q_{\rm max}$,
where $\langle \cdot,e \rangle$ takes the
value $t$.  As $t$ decreases towards
$t_{\rm min} :=\min_\Gamma \langle \cdot,e \rangle$,
the cardinality of the
fibers of $\langle \cdot, e \rangle$ must change from two
to at least three, since otherwise
$\Gamma$ would be homeomorphic to $S^1$.  If the cardinality
increases at another critical point $q_{\rm crit}\neq q^\pm$,
then $q_{\rm crit}$ is a local maximum point and
$\mu(e) \geq
{\rm nlm}(e,q_{\rm max}) + {\rm nlm}(e,q_{\rm crit}) = 2 > 3/2$,
as we wished to show.  The remaining possibility is
that the cardinality increases at a vertex, say $q^+$.  Since the
cardinality increases as $t$ decreases through
$\langle q^+,e \rangle$, we have the down-valence $d^-(e,q^+)$
strictly greater than the up-valence $d^+(e,q^+)$.  By the
Combinatorial Lemma \ref{combin}, ${\rm nlm}(e,q^+)\geq 1/2$, and
in this case also, we have
$\mu(e) \geq
{\rm nlm}(e,q_{\rm max}) + {\rm nlm}(e,q^+) \geq 3/2.$
This shows that for any theta graph $\Gamma$, and almost
any $e \in S^2,$ $\mu(e) \geq 3/2$;  hence by Theorem
\ref{muthm}, the net total curvature $\nc(\Gamma) \geq 3\pi.$

Before discussing the equality case $\nc(\Gamma) = 3\pi$,
we consider the intersection of knots with planes.

%
%
\begin{lem}\label{jordan}
Let $\Gamma_0 \subset \Re^3$ be homeomorphic to $S^1$.  Then
{\bf either} (i) there is a three-parameter family of planes
$P_t^e=\{x\in \Re^3: \langle e,x \rangle = t\}$, 
$t_0<t<t_0+\delta$, $e$ in a neighborhood of 
$e_1 \in S^2$, such that each $P_t^e$ meets
$\Gamma_0$ in at least four points; {\bf or} (ii) $\Gamma_0$ is a
convex plane curve. 
\end{lem}
\pf
If $\Gamma_0$ is not planar, then there exist four non-coplanar
points $p_1, p_2, p_3, p_4$, numbered in order around $\Gamma_0$.
Let an oriented plane $P_0$ be chosen to
contain $p_1$ and $p_3$ and rotated until both $p_2$ and $p_4$ are
above $P_0$ strictly (note that no three of the points can be
collinear).  Write $e_1$ for the unit normal vector to $P_0$ on
the side where $\Gamma_0$ lies.  Then the set $P_t \cap \Gamma_0$
contains at least four points, for $t_0=0 < t < \delta \ll 1$,  
since each plane $P_t=P_t^{e_1}$ meets each of the four open 
arcs between the points $p_1, p_2, p_3, p_4$.
Conclusion (i) remains true, with some $0<t_0 \ll \delta$,  
when the normal vector
$e_1$ to $P_0$ is replaced by any nearby $e \in S^2$.

If $\Gamma_0$ is planar but nonconvex, then there exists a plane 
$P_0=P_0^{e_1}$, transverse to the plane containing $\Gamma_0$, which
supports $\Gamma_0$ and touches $\Gamma_0$ at two distinct points,
but does not include the arc of $\Gamma_0$ between these two points.  
Consider disjoint open arcs of $\Gamma_0$ on either side of these
two points and including points not in $P_0$.  Then for 
$0 < t < \delta \ll 1$, the set $P_t \cap \Gamma_0$ contains at 
least four points, since the planes $P_t=P_t^{e_1}$ meet each 
of the four disjoint arcs.  Here once again $e_1$ may be 
replaced by any nearby unit vector $e$,  and the plane $P_t^e$
will meet $\Gamma_0$ in at least four points, for $t$ in a 
nonempty open interval $t_0<t<t_0+\delta$.
\qed\\

We show next that if a graph $\Gamma$ satisfies alternative (i) in
the conclusion of Lemma \ref{jordan}, then $\mu(e) \geq 2$ for $e$
in a nonempty open set of $S^2$.  Namely, since for
$t_0<t<t_0+\delta$, $P_t^e$ meets $\Gamma$ in at least four points,
the sum of $d^-(e,q)-d^+(e,q)$ over all vertices and critical
points $q$ of $\Gamma$ with 
$\langle e,q \rangle \geq t_0 + \delta$ is at least $4$.  Thus by
Lemma \ref{combin} and Definition \ref{defmu}, $\mu(e) \geq 2$.

Now consider the {\bf equality} case of a theta graph $\Gamma$
with $\nc(\Gamma) = 3\pi$.  As we have seen, the multiplicity 
$\mu(e) = 3/2$ for a.a. $e \in S^2$, so alternative (i) in the 
conclusion of Lemma \ref{jordan} is impossible for $\Gamma$
or for any subgraph $\Gamma_0$ homeomorphic to $S^1$.  By Lemma
\ref{jordan}, such subgraphs $\Gamma_0$ must be planar and convex.

$\Gamma$ consists of three arcs $a_1$, $a_2$ and
$a_3$, meeting at the two vertices $q^+$ and $q^-$. 
The three Jordan curves $\Gamma_1 := a_2 \cup a_3$,
$\Gamma_2 :=  a_3 \cup a_1$ and  $\Gamma_3 := a_2 \cup a_3$
are each planar and convex.  It
follows that $\Gamma_1, \Gamma_2$ and $\Gamma_3$ lie in a common
plane.  In terms of the topology of this plane, one of the three
arcs $a_1$, $a_2$ and
$a_3$ lies in the middle.  But the middle arc must be a line
segment, as it needs to
be a shared piece of two convex curves bounding disjoint open sets
in the plane.  The conclusion is that $\Gamma$ 
is a planar, convex Jordan curve, plus a chord.  This concludes the 
argument in the case of equality $\nc(\Gamma)= 3\pi.$

We finally turn our attention to the F\'ary-Milnor type {\bf upper
bound}, to ensure that a $\theta$-graph is isotopically standard:
we shall assume that $\nc(\Gamma) < 4\pi$.  By Theorem \ref{muthm},
since $S^2$ has area $4\pi$, it follows that there exists
$e_0\in S^2$ with $\mu(e_0) < 2$.  Since $\mu(e_0)$ is a
half-integer, and since $\mu(e) \geq 3/2$,
as we have shown in the first part of this proof,
we have $\mu(e_0) = 3/2$ exactly.

We shall show that there are at most two points  $q \in \Gamma$,
either vertices or critical points, with ${\rm nlm}(e_0,q) >0$.
The ``height" function $\langle \cdot, e_0 \rangle$ has a
maximum at
a point $q_{\rm max}$, with ${\rm nlm}(e_0, q_{\rm max})=1$ or
$3/2$ according as the valence $d(q_{\rm max}) =2$ or $3$.  
If there is a second point $q$ with ${\rm nlm}(e_0,q)>0$, then
${\rm nlm}(e_0,q_{\rm max})=1$ and
${\rm nlm}(e_0, q)=\frac12$, which implies by Lemma \ref{combin} 
that $q$ is a vertex, say $q^+$, with
$d^+(q^+) = 1$ and $d^-(q^+) = 2$ (recall that
the valence $d(q^\pm) = 3$).

       Observe that for any $e\in S^2$,
$\sum_{q}{\rm nlm}(e,q) = 0$;  in fact the number of
local maxima along any parameterization $\Gamma'$ of 
$\widetilde\Gamma$ is equal to the number of local minima.  
But $[{\rm nlm}(-e_0,q) = -[{\rm nlm}(e_0,q)$.  It follows 
that $\mu(-e_0) = \sum [{\rm nlm}(e_0,q)]^-  = 3/2$.
We may apply the argument above, replacing $e_0$ with $-e_0$, to
show that either ${\rm nlm}(e_0, q_{\rm min})= -3/2$ and
${\rm nlm}(e_0, q) \geq 0$ elsewhere; or
${\rm nlm}(e_0, q_{\rm min})= -1$ and and there is a second point,
which must be a vertex, say $q^-$, with
${\rm nlm}(e_0, q^-) = -\frac12.$

Thus, according to Lemma \ref{combin}, there are four cases,
depending on the valences $d(q_{\rm max})$ and
$d(q_{\rm min}) \in \{2,3\},$
for the cardinality $\#(e_0,t)$ of the fiber
$\{q\in \Gamma: \langle e_0,q \rangle =t \}$ as $t$ decreases from
$t_{\rm max}$  to $t_{\rm min}$.  Write
$t^\pm := \langle q^\pm, e_0 \rangle$.  The four cases, listed by
$\left(d(q_{\rm max}),d(q_{\rm min})\right)$, are:
\begin{itemize}
\item[(3,3):]  $\#(e_0,t) \equiv 3$,
  $t_{\rm min} = t^- < t <t_{\rm max}=t^+$.
\item[(3,2):]  $\#(e_0,t) = 3$ for
  $t^- <t<t_{\rm max}= t^+;  \#(e_0,t)=2$ for $t_{\rm min} < t < t^-$.
\item[(2,3):]  $\#(e_0,t) =2$ for
  $t^+ < t <t_{\rm max}$; $\#(e_0,t) = 3$ for 
$t_{\rm min}=t^- < t < t^+$.
\item[(2,2):]  $\#(e_0,t) =2$ for
  $t^+<t<t_{\rm max}$ and for  $t_{\rm min} < t < t^-$;
  $\#(e_0,t) = 3$ for $t^- < t < t^+$.
\end{itemize}

In each of these four cases,
we shall show that $\Gamma$ is isotopic in $\Re^3$ to a planar
$\theta$-graph via an isotopy which does not change the values of
$\langle e_0, \cdot \rangle$.  Let
us consider the fourth case {\bf $(2,2)$} in detail, and observe
that the other three cases follow in a similar fashion.

Write $P_t=P_t^{e_0}$ for the ``horizontal" plane of $\Re^3$ at
height $t$.  Then $\#(e_0,t)$ is the cardinality of $P_t \cap \Gamma$.

In the fourth case {\bf $(2,2)$}, for $t^-<t<t^+$, there are
$\#(e_0,t) = 3$ points in $\Gamma \cap P_t$.  One of the three
edges of \, $\Gamma$, say $a_2$,
lies entirely in the closed slab between $P_{t^-}$ and $P_{t^+}$.
As $t \rightarrow t^\pm$, two of the points converge
from a well-defined direction, since the three unit tangent
vectors at $q^\pm$ are distinct.  Letting $t$ decrease from
$t=t^+$, there is an isotopy of the plane $P_t$, varying
continuously with $t$, so that the three points of
$\Gamma \cap P_t$ become collinear, with the points
$a_1 \cap P_t$, $a_2 \cap P_t$ and $a_3 \cap P_t$ appearing
in that order, or the reverse order, along a line.
After a further isotopy in $\Re^3$ which translates and rotates
each plane $P_t$ rigidly, we may achieve that
$P_t \cap \Gamma \subset Q$ for some plane $Q$ in $\Re^3$
transverse to the planes $P_t$.

For $t>t^+$ and for $t< t^-$,
we may continuously rotate and translate the planes $P_t$
so that this ``top" portion and this ``bottom" portion of
$\Gamma$ each lie in the same plane $Q$.
Thus, $\Gamma$ is isotopic to a planar graph.

We may now find a further isotopy of the plane $Q$, extendible to
$\Re^3$, which deforms $\Gamma$
into a circle with one diameter.  Therefore, $\Gamma$ is isotopic 
in $\Re^3$ to the standard $\theta$ graph.
\qed\\

It might be noted that the last part of our proof does
{\em not} hold when $\Gamma$ has the homeomorphism type of
a circle with {\em two} disjoint, parallel
chords.  Of course, the Fenchel-type lower bound for $\nc(\Gamma)$
is larger than $3\pi$, namely $4\pi$.  However, this graph may be
embedded into $\Re^3$ so that it is not isotopic to a planar
graph, but so that a particular vector $e_0$ has the
minimum value $\mu(e_0) = 2$:  think of twisting the two chords of
the standard graph an even number of times about each other, 
without introducing any
critical points of $\langle e_0,t \rangle$, and leaving the
rest of the graph alone.  Thus, one may construct nonisotopic
graphs with $\nc(\Gamma)$ slightly greater than $4\pi$.  
Nonetheless, we conjecture that any
graph $\Gamma$ of this homeomorphism type with $\nc(\Gamma) =
4\pi$ is a convex plane curve with two disjoint chords.

There remains a challenging question to determine the
F\'ary-Milnor type lower bound for the net total curvature of
graphs in $\Re^3$ of {\it any} specific homeomorphism type which
are not isotopic to a standard embedding of that graph.

%
\section{Nowhere-smooth graphs}\label{nonsmooth}

Milnor extended his proof of the F\'ary-Milnor Theorem to
continuous knots in \cite{M};  we shall carry out an analogous
extension to continuous graphs.

Up to this point, we have treated graphs which are
piecewise $C^2$, so that the definition of net
total curvature in Definition \eqref{defnet} makes sense:
$$
{\nc}(\Gamma):=
\sum_{i=1}^N {\rm nc}(q_i) + \int_{\Gamma_{\rm reg}} |\vec{k}| \, ds
$$
Recall that Milnor~\cite{M} defines the total
curvature of a continuous simple closed curve $C$ as the supremum
of the total curvature of all polygons inscribed in $C$.  In
analogy to Milnor, we define net total curvature of a 
{\it continuous} graph $\Gamma$ to be the supremum of the net total
curvature of all polygonal graphs $P$ suitably inscribed in
$\Gamma$ as follows.  

%
%
\begin{defi}\label{gamapprox}
For a given continuous graph $\Gamma$, we
say a polygonal graph $P \subset \Re^3$ is 
{\em $\Gamma$-approximating}, provided that its topological
vertices (those of valence $\neq 2$) are exactly the topological
vertices of \, $\Gamma$, and having the same valences; and that an 
arc of $P$ between two topological vertices corresponds one-to-one 
to each edge of \, $\Gamma$ between those two vertices, the vertices of
valence $2$ of that
arc lying in order along the corresponding edge of \, $\Gamma$. 
\end{defi}
Note that if $P$ is a $\Gamma$-approximating polygonal graph,
then $P$ is homeomorphic to $\Gamma$.   Recall that
according to Proposition \ref{monotmu}, if $P$ and 
$\widetilde{P}$ are $\Gamma$-approximating polygonal
graphs, and $\widetilde{P}$ is a refinement of $P$, then 
$\nc(\widetilde{P}) \geq \nc(P)$.

%
\begin{defi}\label{gendefnet}
Define the {\em net total curvature} of a continuous graph 
$\Gamma$ by
\[
\nc(\Gamma) := \sup_{P} \nc(P)
\]
where the supremum is taken over all $\Gamma$-approximating
polygonal graphs $P$.  
\end{defi}

Here $\nc(P)$ is given as in definition \ref{defnet}, the 
definition of net total curvature, by
$$
{\nc}(P):=
\sum_{i=1}^N {\rm nc}(q_i), 
$$
where $q_1, \dots, q_N$ are the vertices of $P$.

Definition \ref{gendefnet} is consistent with Definition
\ref{defnet} in the case where both may be applied, namely in the
case of a piecewise $C^2$ graph $\Gamma$.  Namely, Milnor 
showed that the total curvature ${\mathcal C}(\Gamma_0)$ of 
a smooth curve $\Gamma_0$ is 
the supremum of the total curvature of inscribed polygons 
(\cite{M}, p. 251), which gives the required supremum for each
edge.  At a vertex $q$ of the piecewise-$C^2$ graph $\Gamma$, 
as a sequence $P_k$ of  $\Gamma$-approximating 
polygons become arbitrarily fine, a vertex $q$ of $P_k$
(and of $\Gamma$) has unit tangent vectors converging in $S^2$
to the unit tangent vectors to $\Gamma$ at $q$.  
It follows that
for $1\leq i\leq d(q)$, $\chi_i^{P_k} \to \chi_i^{\Gamma}$ in
measure on $S^2$, and therefore
${\rm nc}_{P_k}(q_k) \to {\rm nc}_\Gamma(q)$.

%
\begin{defi}\label{critpt}
We say a point $q\in \Gamma$ is {\em critical} 
relative to $e \in S^2$ when 
$q$ is a topological vertex of \,$\Gamma$ or when
$\langle e, \cdot\rangle$ is not monotone in any open 
interval of \, $\Gamma$ containing $q$.
\end{defi}

Note that at some points of a differentiable curve,
$\langle e, \cdot\rangle$ may have derivative zero but still
not be considered a critical point relative to $e$ by our 
definition.  This is appropriate to the $C^0$ category.
For a continuous graph $\Gamma$, when 
$\nc (\Gamma)$ is finite, we shall 
show that the number of critical points is finite for almost all
$e$ in $S^2$ (see Lemma \ref{fincrit} below).   

%
\begin{lem}\label{monotcvge}
Let $\Gamma$ be a continuous, finite graph in $\Re^3$, and choose
a sequence $\widehat{P_k}$ of \, $\Gamma$-approximating polygonal 
graphs with $\nc(\Gamma)= \lim_{k \rightarrow \infty}
\nc(\widehat{P_k}).$ 
Then for each
$e \in S^2$, there is a refinement $P_k$ of $\widehat{P_k}$ 
such that $\lim_{k \rightarrow \infty}\mu_{P_k}(e)$ exists in
$[0,\infty]$. 
\end{lem}
\pf 
As a first step, for each $k$ in sequence, we refine 
$\widehat{P_k}$ to include all vertices of
$\widehat{P_{k-1}}$.  Then for all $e\in S^2$, 
$\mu_{\widehat{P_k}}(e) \geq \mu_{\widehat{P_{k-1}}}(e)$,
by Proposition \ref{monotmu}.
As the second step, we refine $\widehat{P_k}$ so that the arc of
$\Gamma$ corresponding to each edge of $\widehat{P_k}$ has
diameter $\leq 1/k$.  As the third step, given a particular 
$e \in S^2$, for each edge $\widehat{E_k}$ of $\widehat{P_k}$, 
we add $0,1$ or $2$ points from $\Gamma$ as vertices of 
$\widehat{P_k}$ so that 
$\max_{\widehat{E_k}}\langle e,\cdot\rangle=
\max_E\langle e,\cdot\rangle$
where $E$ is the closed arc of \, $\Gamma$ corresponding to
$\widehat{E_k}$;  
and similarly so that 
$\min_{\widehat{E_k}} \langle e,\cdot\rangle=
\min_E\langle e,\cdot\rangle$.
Write $P_k$ for the result of this three-step refinement.  
Note that all vertices of $P_{k-1}$ appear among the vertices
of $P_k$.  Then by Proposition \ref{monotmu}, 
$$ \nc(\widehat{P_k}) \leq \nc(P_k) \leq \nc(\Gamma), $$
so we still have 
$\nc(\Gamma)= \lim_{k \rightarrow \infty} \nc(P_k).$

Now compare the values of
$\mu_{P_k}(e)=\sum_{q\in P_k} {\rm nlm_{P_k}}^+(e,q)$ with
the same sum for $P_{k-1}$.
Since $P_k$ is a refinement of $P_{k-1}$, we have
$\mu_{P_k}(e) \geq \mu_{P_{k-1}}(e)$ by Proposition
\ref{monotmu}. 

Therefore the values $\mu_{P_k}(e)$ are non-decreasing in $k$, which
implies they are either convergent or properly divergent;  in the
latter case we write 
$\lim_{k \rightarrow \infty}\mu_{P_k}(e)= \infty$.
\qed

%
\begin{defi}
For a continuous graph $\Gamma$, define the {\em multiplicity} at
$e\in S^2$ as $\mu_\Gamma(e):= 
\lim_{k \rightarrow \infty}\mu_{P_k}(e) \in [0,\infty]$, 
where $P_k$ is a sequence of  $\Gamma$-approximating 
polygonal graphs, refined with respect to $e$, as given 
in Lemma \ref{monotcvge}.
\end{defi}

%
\begin{rem}
Note  that any two $\Gamma$-approximating polygonal graphs 
have a common refinement.  Hence, from the proof of Lemma 
\ref{monotcvge}, any two choices of sequences 
$\{\widehat{P_k}\}$ of $\Gamma$-approximating
polygonal graphs lead to the same value $\mu_\Gamma(e)$.
\end{rem}

%
\begin{lem}\label{a.a.e}
Let $\Gamma$ be a continuous, finite graph in $\Re^3$. 
Then $\mu_\Gamma: S^2 \to [0,\infty]$ takes its values in the
half-integers, or $+ \infty$.
Now assume $\nc(\Gamma) < \infty$.  Then $\mu_\Gamma$ is integrable, 
hence finite almost everywhere on $S^2$, and 
\begin{equation}\label{nc=int}
 \nc(\Gamma) = \frac12 \int_{S^2} \mu_\Gamma(e) \, dA_{S^2}(e).
\end{equation}
For almost all $e \in S^2$, a sequence $P_k$ of
\, $\Gamma$-approximating polygonal graphs may be chosen 
(depending on $e$) so that each
local extreme point $q$ of $\langle e,\cdot\rangle$ along $\Gamma$ 
occurs as a vertex of $P_k$ for sufficiently large $k$.
\end{lem}
\pf
The half-integer-valued functions $\mu_{P_k}$ are non-negative,
integrable on $S^2$ with bounded integrals since 
$\nc(\Gamma) < \infty$, and monotone increasing in $k$.  
Thus for almost all $e \in S^2$,
$\mu_{P_k}(e) = \mu_\Gamma(e)$ for $k$ sufficiently large.
It follows that if $\mu_\Gamma(e)$ is finite, it must be a 
half-integer.

Since the functions $\mu_{P_k}$ are non-negative and pointwise
non-decreasing almost everywhere on $S^2$, it now follows 
from the Monotone Convergence Theorem that 
$$\int_{S^2} \mu_\Gamma(e)\, dA_{S^2}(e) = 
\lim_{k\to \infty}\int_{S^2} \mu_{P_k}(e)\, dA_{S^2}(e)=
2 \nc(\Gamma).$$

Finally, let the $\Gamma$-approximating polygonal graphs
$P_k$ be chosen to have maximum edge length $\to 0$.  
For almost all $e \in S^2$, $\langle e,\cdot\rangle$ is
not constant along any open arc of  $\Gamma$, and 
$\mu_\Gamma(e)$ is finite.  Given such an $e$,
choose $\ell = \ell(e)$ sufficiently 
large that $\mu_{P_k}(e) = \mu_\Gamma(e)$ and 
$\mu_{P_k}(-e) = \mu_\Gamma(-e)$ for all $k \geq \ell$.  Then for
$k \geq \ell$, along any edge $E_k$ of $P_k$ with corresponding
arc $E$ of $\Gamma$, the maximum and
minimum values of $\langle e,\cdot\rangle$ along $E$ occur at the
endpoints, which are also the endpoints of $E_k$.  Otherwise, as
$P_k$ is further refined, new interior local maximum resp. local
minimium points of $E$ would
contribute a new, positive value to $\mu_{P_k}(e)$ resp. to
$\mu_{P_k}(-e)$ as $k$
increases.  Since the diameter of the corresponding arc $E$ of
$\Gamma$ tends to zero as $k \to \infty$, any local
maximum or local minimum of $\langle e,\cdot\rangle$ must become
an endpoint of some edge of $P_k$ for $k$ sufficiently large, and
for $k\geq \ell$ in particular.
\qed\\

Our last lemma focuses on the regularity of a graph $\Gamma$,
originally only assumed continuous, provided it has finite net
total curvature, or another notion of total curvature of a graph 
which includes the total curvature of the edges.

%
\begin{lem}\label{fincrit}
Let $\Gamma$ be a continuous, finite graph in $\Re^3$, with
$\nc(\Gamma)<\infty$. Then $\Gamma$ has continuous one-sided unit
tangent vectors $T_\pm(p)$ at each point $p$.  If $p$ is a
vertex of valence $d$, then each of the $d$ edges which meet at
$p$ have well-defined unit tangent vectors at $p$:
$T_1(p),\dots,T_d(p)$.  For almost all $e \in S^2$, 
\begin{equation}\label{mu=sum}
\mu_\Gamma(e) = \sum_q\{{\rm nlm}(e,q)\}^+,
\end{equation}
where the sum is over the {\em finite} number of topological
vertices of \, $\Gamma$ and critical points $q$ of 
$\langle e, \cdot \rangle$ along $\Gamma$. 
Further, for each $q$,
${\rm nlm}(e,q)= \frac12[d^-(e,q) - d^+(e,q)]$.
All of these critical points which are not topological vertices 
are local extrema of $\langle e,\cdot\rangle$ along $\Gamma$.  
\end{lem}
\pf
We have seen in the proof of Lemma \ref{a.a.e} that for almost all 
$e \in S^2$, the linear function $\langle e,\cdot\rangle$ is 
not constant along any open arc of  $\Gamma$, and
there is a sequence $\{P_k\}$ of $\Gamma$-approximating 
polygonal graphs with
$\mu_\Gamma(e) = \mu_{P_k}(e)$ for $k$ sufficiently large.
We have further shown that each local maximum point of 
$\langle e,\cdot\rangle$ is a vertex of $P_k$ for $k$ large
enough.  Recall that
$\mu_{P_k}(e) = \sum_q{\rm nlm}_{P_k}^+(e,q)$.
Thus, each local maximum point $q$ for 
$\langle e, \cdot \rangle$ along $\Gamma$ provides a non-negative
term ${\rm nlm}^+(e,q)$ in the sum for $\mu_{P_k}$.  Fix such an
integer $k$.

Consider a point $q\in \Gamma$ which is not a topological 
vertex of $\Gamma$ but is a critical point of 
$\langle e,\cdot\rangle$.  We shall show, by an argument similar 
to one used by van Rooij in \cite{vR}, that $q$ 
must be a local extreme point.
As a first step, we show that $\langle e,\cdot\rangle$ is monotone
on a sufficiently small interval on either side of $q$.
Choose an ordering of the closed edge $E$ of 
$\Gamma$ containing $q$, and consider the interval $E_+$ of 
points $\geq q$ with respect to this ordering.  Suppose that 
$\langle e,\cdot\rangle$ is not monotone on any subinterval of
$E_+$ with $q$ as endpoint.
Then in any interval $(q,r_1)$ there are points $p_2 > q_2 > r_2$
so that the numbers 
$\langle e,p_2 \rangle,\langle e,q_2\rangle,\langle e, r_2\rangle$
are not monotone.
It follows by an induction argument that there exist decreasing 
sequences $p_n \to q$, $q_n \to q$, and 
$r_n \to q$ of points of $E_+$ such that for each $n$, 
$r_{n-1} > p_n > q_n > r_n > q$, but the value 
$\langle e,q_n\rangle$ lies outside of the closed interval between 
$\langle e,p_n\rangle$ and $\langle e,r_n\rangle$.  As a consequence, 
there is a local extremum $s_n \in (r_n, p_n)$.  Since 
$r_{n-1} > p_n$, the $s_n$ are all distinct, $1\leq n < \infty$.  
But by Lemma \ref{a.a.e},
all local extreme points, specifically $s_n$, of 
$\langle e, \cdot \rangle$ along $\Gamma$ 
occur among the {\em finite}
number of vertices of $P_k$, a contradiction.  This shows that 
$\langle e, \cdot \rangle$ is monotone on an interval to the right
of $q$.  A similar argument shows that $\langle e, \cdot \rangle$
is monotone on an interval to the left of $q$.

Recall that  for a {\em critical point} $q$ relative to $e$, 
$\langle e,\cdot\rangle$ is not monotone on any neighborhood 
of $q$.  
By this definition, the sense of monotonicity of 
$\langle e, \cdot \rangle$ must be opposite on the two
sides of $q$.  Therefore every critical point $q$, which is not a
topological vertex, is a local extremum.

Now choose $k$ large enough that $\mu_\Gamma(e) = \mu_{P_k}(e)$.  Then for
any edge $E_k$ of $P_k$, the function $\langle e, \cdot \rangle$
is monotone along the corresponding arc $E$ of \, $\Gamma$, as well
as along $E_k$.  Also, $E$ and $E_k$ have common end points.
It follows that for each $t \in \Re$, the cardinality  
$\#(e,t)$ of the fiber $\{q\in \Gamma: \langle e,q \rangle =t \}$ 
is the same for $P_k$ as for $\Gamma$.  We may see from 
Lemma \ref{combin} applied to $P_k$ that for each vertex 
or critical point $q$,
${\rm nlm}(e,q) = \frac12[d^-(e,q) - d^+(e,q)]$;  but 
${\rm nlm}(e,q)$ and $d^\pm(e,q)$ have the {\it same} values 
for $\Gamma$ as for $P_k$.  Finally, the formula 
$\mu_\Gamma(e) = \sum_q\{{\rm nlm}_\Gamma(e,q)\}^+$ now follows 
from the analogous formula for $P_k$, for almost all $e \in S^2$.
\qed

%
\begin{cor}\label{ctsnc}
Let $\Gamma$ be a continuous, finite graph in $\Re^3$, with
$\nc(\Gamma)<\infty$. Then for each  point $q$ of
$\Gamma$, the contribution at $q$ to net total curvature is
given by equation \eqref{defnc}, where for $e \in S^2$, 
$\chi_i(e)=$ the sign of $\langle -T_i(q), e \rangle$, 
$1\leq i \leq d(q)$.  
(Here, if $q$ is not a topological vertex, 
we understand $i \in \{+,-\}$.)
\end{cor}
\pf
According to Lemma \ref{fincrit}, for $1\leq i \leq d(q)$,
$T_i(q)$ is defined and tangent to an edge $E_i$ of $\Gamma$,
which is continuously differentiable at its end point $q$.  
If $P_n$ is a sequence of
$\Gamma$-approximating polygonal graphs with maximum edge length
tending to $0$, then the corresponding unit tangent vectors
$T^{P_n}_i(q) \to T^{\Gamma}_i(q)$ as $n \to \infty$.  For each
$P_n$, we have 
$$
{\rm nc}^{P_n}(q) =
\frac12\int_{S^2}
\left[\sum_{i=1}^d{\chi_i}^{P_n}(e)\right]^+\,dA_{S^2}(e),
$$
and ${\chi_i}^{P_n} \to {\chi_i}^\Gamma$ in measure on $S^2$.
Hence, the integrals for $P_n$ 
converge to those for $\Gamma$, which is equation \eqref{defnc}.
\qed\\

We are ready to state the formula for net total curvature,
by localization on $S^2$, a generalization of Theorem \ref{muthm}:

%
%
\begin{thm}\label{muthm2}
For a continuous graph $\Gamma,$ the net total curvature 
$\nc(\Gamma) \in (0,\infty]$ has 
the following representation:
$$ \nc(\Gamma) = \frac12 \int_{S^2} \mu(e) \,dA_{S^2}(e), $$
where, for almost all $e \in S^2$, the multiplicity 
$\mu(e)$ is a positive half-integer or $+\infty$, given as the finite
sum \eqref{mu=sum}. 
\end{thm}

\pf
If $\nc(\Gamma)$ is finite, then the theorem follows from Lemma 
\ref{a.a.e} and Lemma \ref{fincrit}.

Suppose $\nc(\Gamma) = \sup \nc(P_k)$ is infinite, where $P_k$ 
is a refined sequence of polygonal graphs as in Lemma
\ref{monotcvge}.  Then $\mu_\Gamma(e)$ is the non-decreasing 
limit of $\mu_{P_k}(e)$ for all $e \in S^2$.  Thus    
$\mu_\Gamma(e) \geq \mu_{P_k}(e)$ for all $e$ and $k$, and
$\mu_\Gamma(e) = \mu_{P_k}(e)$ for $k\geq\ell(e)$.  This implies
that $\mu_\Gamma(e)$ is a positive half-integer or $\infty$.
Since $\nc(\Gamma) = \infty$, the integral 
$$ \nc(P_k) = \frac12 \int_{S^2} \mu_{P_k}(e) \,dA_{S^2}(e)$$
is arbitrarily large as $k \to \infty$, but for each $k$ is 
less than or equal to 
$$ \frac12 \int_{S^2} \mu_\Gamma(e) \,dA_{S^2}(e).$$
Therefore this latter integral equals $\infty$, and thus equals
$\nc(\Gamma).$
\qed\\

We turn our attention next to the tameness of graphs of finite
total curvature.

%
\begin{pro}\label{untangle}
Let $n$ be a positive integer, and write $Z$ for the set of $n$th
roots of unity in $\C = \Re^2$.  Given a continuous one-parameter 
family $S_t$, $ 0 \leq t < 1$, of sets of $n$ points in 
$\Re^2$, there exists a continuous
one-parameter family $\Phi_t:\Re^2 \to \Re^2$ of homeomorphisms
with compact support such that $\Phi_t(S_t) = Z$, $0 \leq t < 1$. 
\end{pro}

\pf
It is well known that there is an isotopy $\Phi_0:\Re^2 \to \Re^2$
such that $\Phi_0(S_0) = Z$ and $\Phi_0 =$ id outside of a
compact set.  Write $B_\varepsilon (Z)$ for the union of balls
$B_\varepsilon (\zeta_i)$ centered at the $n$ roots of unity
$\zeta_1, \dots \zeta_n$.  For $\varepsilon < \sin{\frac{\pi}{n}},$
these balls are disjoint.  This completes the case 
$ t_0 = 0$ of the following continuous induction argument.

Suppose that $[0,t_0] \subset [0,1)$ is a subinterval such that
there exists a continuous one-parameter family
$\Phi_t:\Re^2 \to \Re^2$ of homeomorphisms with compact support,
with $\Phi_t(S_t) = Z$ for all $0 \leq t \leq t_0$.  
We shall extend this property to an interval $[0,t_0+\delta]$.
We may choose $0 < \delta < 1-t_0$  such that
$\Phi_{t_0}(S_t) \subset B_\varepsilon(Z)$ for all
$t_0 \leq t \leq t_0 + \delta.$  Write the points of $S_t$ as
$x_i(t), \ 1 \leq i \leq n,$ where
$\Phi_{t_0}(x_i(t)) \in B_\varepsilon(\zeta_i)$.  For each
$t \in [t_0, t_0 + \delta],$ each of the balls
$B_\varepsilon(\zeta_i)$ may be mapped onto itself by a 
homeomorphism $\psi_t$, varying continuously with $t$, such that
$\psi_{t_0}$ is the identity, $\psi_t$ is the identity near the
boundary of $B_\varepsilon(\zeta_i)$ for all 
$t \in [t_0, t_0 + \delta]$, and
$\psi_t(\Phi_{t_0}(x_i(t))) = \zeta_i$ for all such $t$. 
For example, we may construct $\psi_t$ so that
for each $y \in B_\varepsilon(\zeta_i)$, $y-\psi_t(y)$ is
parallel to $\Phi_{t_0}(x_i(t)) - \zeta_i$.  We now define
$\Phi_t = \psi_t \circ \Phi_{t_0}$ for each
$t \in [t_0, t_0 + \delta].$

As a consequence, we see that there is no maximal interval
$[0,t_0] \subset [0, 1)$ such that
there is a continuous one-parameter family
$\Phi_t:\Re^2 \to \Re^2$ of homeomorphisms with compact support
with $\Phi_t(S_t) = Z$, for all $0 \leq t \leq t_0$.
Thus, this property holds for the entire interval $0 \leq t<1$.
\qed\\

In the following theorem, the total curvature of a graph may be
understood in terms of any definition which includes the total
curvature of edges and which is continuous as a function of the
unit tangent vectors at each vertex.  This includes net total
curvature,  maximal total curvature (TC of Taniyama, \cite{T}) and
the cone total curvature of \cite{GY}.

%
\begin{thm}\label{tame}
Suppose $\Gamma \subset \Re^3$ is a continuous graph with finite
total curvature.  Then for any $\varepsilon > 0$, $\Gamma$ is 
isotopic to a $\Gamma$-approximating polygonal graph $P$ with
edges of length at most $\varepsilon$, whose total curvature is 
less than or equal to that of $\Gamma$.
\end{thm}
\pf
Since $\Gamma$ has finite total curvature, by Lemma \ref{fincrit},
at each topological vertex of valence $d$ the edges have
well-defined unit tangent vectors $T_1, \dots, T_d$, which are
each the limit of the unit tangent vectors to the corresponding edges.
If at each vertex the unit tangent vectors $T_1, \dots, T_d$ 
are distinct, then any sufficiently fine
$\Gamma$-approximating polygonal graph will be isotopic to
$\Gamma$;  this easier case is proven.

We consider therefore $n$ edges $E_1, \dots, E_n$ which end at a
vertex $q$ with common unit tangent vectors $T_1 = \dots = T_n$.
Choose orthogonal coordinates $(x,y,z)$ for $\Re^3$ so that this 
common tangent vector $T_1 = (0,0,-1)$ and $q = (0,0,1)$.
After rescaling about $q$ by a sufficiently large factor 
$> \frac{1}{\varepsilon}$, $E_1, \dots, E_n$ form a 
braid $B$ of $n$ strands in the 
slab $0 \leq z < 1$ of $\Re^3$, plus the point $q=(0,0,1)$.  Each
strand $E_i$  has $q$ as an endpoint, and the coordinate $z$ is 
strictly monotone along $E_i$, $1\leq i \leq n$.  Write 
$S_t = B \cap \{ z = t\}$.
Then $S_t$ is a set of $n$ distinct points in the plane 
$\{ z = t\}$ for each $0\leq t<1$.  According to Proposition
\ref{untangle}, there are homeomorphisms $\Phi_t$ of the plane
$\{ z = t\}$ for each $0\leq t<1$, isotopic to the identity in
that plane,  continuous as a function of $t$,
such that $\Phi_t(S_t) = Z \times \{t\},$ where $Z$ is the set of
$n$th roots of unity in the $(x,y)$-plane, and $\Phi_t$ is the
identity outside of a compact set of the plane $\{ z = t\}$.  

We may suppose that $S_t$ lies in the open disk of radius $a(1-t)$
of the plane $\{ z = t\}$, for some (arbitrarily small) 
constant $a>0$. We modify
$\Phi_t$, first replacing its values with $(1-t) \Phi_t$ inside
the disk of radius $a(1-t)$.  We then modify $\Phi_t$ outside the
disk of radius $a(1-t)$, such that $\Phi_t$ is the identity
outside the disk of radius $2a(1-t)$.

Having thus modified the homeomorphisms $\Phi_t$ of the planes
$\{ z = t\}$, we may now define an isotopy $\Phi$ of $\Re^3$ by
mapping each plane $\{ z = t\}$ to itself by the homeomorphism
$\Phi_0^{-1} \circ \Phi_t$, $0\leq t<1$; and extend to the
remaining planes $\{ z = t\}$, $t\geq 1$ and $t<0$, by the 
identity.  Then the closure
of the image of the braid $B$ is the union of line segments from
$q =(0,0,1)$ to the $n$ points of $S_0$ in the plane $\{ z = 0\}$.
Since each $\Phi_t$ is isotopic to the identity
in the plane $\{ z = t\}$, $\Phi$ is isotopic to
the identity of $\Re^3$.

This procedure may be carried out in disjoint sets of $\Re^3$
surrounding each unit vector which occurs as tangent 
vector to more than one edge at a vertex of $\Gamma$.  Outside
these sets, we inscribe a polygonal arc in each edge of $\Gamma$
to obtain a $\Gamma$-approximating polygonal graph $P$.  By
Definition \ref{gendefnet}, $P$ has total curvature less than or
equal to the total curvature of $\Gamma$.

\qed\\

Artin and Fox \cite{AF} introduced the notion of {\em tame} and 
{\em wild}
knots in $\Re^3$;  the generalization to graphs is the following

%
\begin{defi}
We say that a graph in $\Re^3$ is {\em tame} if it is isotopic 
to a polyhedral graph;  otherwise, it is {\em wild}.
\end{defi}

Milnor proved in \cite {M} that knots of finite total curvature
are tame. More generally, we have

%
\begin{cor}
A continuous graph $\Gamma \subset \Re^3$ of finite total
curvature is tame.
\end{cor}
\pf
This is
an immediate consequence of Theorem \ref{tame}, since the
$\Gamma$-approximating polygonal graph $P$ is isotopic to
$\Gamma$.
\qed

%
\begin{thm}\label{theta2}
Suppose $\Gamma \subset \Re^3$ is a continuous theta graph.
Then $\nc(\Gamma) \geq 3\pi$.  If $\nc(\Gamma) < 4\pi$, then
$\Gamma$ is isotopic in $\Re^3$ to the standard theta graph.
Moreover, when $\nc(\Gamma) = 3\pi$, the graph is a planar convex
curve plus a straight chord.
\end{thm}
\pf
It follows from Theorem \ref{tame} that for any $\theta$-graph
$\Gamma$ of finite net total curvature, there is a
$\Gamma$-approximating polygonal $\theta$-graph $P$ isotopic to
$\Gamma$, with $\nc(P) \leq \nc(\Gamma)$ and as close as desired
to $\nc(\Gamma)$.

If a $\theta$-graph 
$\Gamma$ would have $\nc(\Gamma) < 3\pi$, then the $\Gamma$-approximating
polygonal graph $P$ would also have $\nc(P) < 3\pi$, in contradiction to
Theorem \ref{thetathm}.  This shows that $\nc(\Gamma) \geq 3\pi$.

If equality $\nc(\Gamma) = 3\pi$ holds, then 
$\nc(P) \leq \nc(\Gamma) = 3\pi$,
so that by Theorem \ref{thetathm}, $\nc(P)$ must equal $3\pi$, and 
$P$ must be a convex planar curve plus a chord.  But this holds
for {\em all} $\Gamma$-approximating polygonal graphs $P$, implying that 
$\Gamma$ itself must be a convex planar curve plus a chord.

Finally, if $\nc(\Gamma) < 4\pi$, then $\nc(P) < 4\pi$, implying
that $P$ is isotopic to the standard $\theta$-graph.  But $\Gamma$
is isotopic to $P$, and hence is standard.
\qed

\bigskip

\normalsize

\begin{tabbing}
aaaaaaaaaaaaaaaaaaaaaaaaaaaaaasssssssssssssssss \=
bbbbbbbbbbbbbbbbbbbbbbbbbbbbbbb\kill

Robert Gulliver\> Sumio Yamada\\
School of Mathematics \>Mathematical Institute\\
University of Minnesota \> Tohoku University\\
Minneapolis MN 55414 \> Aoba, Sendai, Japan 980-8578\\
{\tt gulliver@math.umn.edu} \> {\tt yamada@math.tohoku.ac.jp}\\
{\tt www.ima.umn.edu/\~{ }gulliver}\> {\tt www.math.tohoku.ac.jp/\~{ }yamada}
\\
\end{tabbing}

\end{document}